\documentclass[a4paper,12pt]{amsart}

\usepackage{amssymb,enumerate,psfrag,graphicx,amsfonts,amsthm,mathrsfs,amsmath,amscd,graphicx}
\usepackage{cleveref}
\usepackage{xcolor}
\usepackage{tikz}
\usetikzlibrary{calc}
\usetikzlibrary{arrows.meta, positioning, calc}
\usetikzlibrary{shapes.geometric}
\usetikzlibrary{arrows}
\usepackage{pgffor}
\usetikzlibrary{decorations.markings}
\usetikzlibrary{arrows.meta, decorations.pathreplacing}
\usetikzlibrary{decorations.text, positioning}
\tikzset{negated/.style={
        decoration={markings,
            mark= at position 0.5 with {
                \node[transform shape] (tempnode) {$\backslash$};
            }
        },
        postaction={decorate}
    }
}
\usetikzlibrary{decorations.markings}
\usetikzlibrary{decorations.text}
\usetikzlibrary{positioning}
\tikzset{
    vertex/.style={draw,circle,thick, inner sep=1.5pt,minimum size=6pt,
    },
    edge/.style={thick},
    dedge/.style = {->,> = latex',thick}
}

\definecolor{mygray}{gray}{0.8}

\newtheorem{theo}{Theorem}[section]
\newtheorem{theorem}[theo]{Theorem}
\newtheorem{corollary}[theo]{Corollary}
\newtheorem{lemma}[theo]{Lemma}
\newtheorem{proposition}[theo]{Proposition}

\newtheorem{theor}{Theorem}

\theoremstyle{definition}

\newtheorem{example}[theo]{Example}
\newtheorem{remark}[theo]{Remark}
	
\newtheorem{problem}[theo]{Problem}	
	
 \newtheorem{definition}[theo]{Definition}

\numberwithin{equation}{section}

\title[Switching reconstruction of signed graphs: net-degree and balance]{Switching reconstruction of signed graphs: net-degree and balance}

\author[M. Cavaleri]{Matteo Cavaleri}
\address{Matteo Cavaleri, Universit\`{a}  Niccol\`{o} Cusano, Dipartimento di Ingegneria - Via Don Carlo Gnocchi, 3 00166 Roma, Italia}
\email{matteo.cavaleri@unicusano.it (Corresponding author)
}

\author[A. Donno]{Alfredo Donno}
\address{Alfredo Donno, Universit\`{a}  Niccol\`{o} Cusano, Dipartimento di Ingegneria - Via Don Carlo Gnocchi, 3 00166 Roma, Italia}
\email{alfredo.donno@unicusano.it}

\author[S. Spessato]{Stefano Spessato}
\address{Stefano Spessato, Universit\`{a}  Niccol\`{o} Cusano, Dipartimento di Ingegneria - Via Don Carlo Gnocchi, 3 00166 Roma, Italia}
\email{stefano.spessato@unicusano.it}

\begin{document}

\maketitle

\begin{abstract}
The switching reconstruction problem for $2$-edge-colored graphs, introduced by
Bondy and Mercier, asks whether such a graph is determined up to isomorphism by
the multiset of the graphs obtained by exchanging the colors of the edges at 
each vertex; for complete underlying graphs it is the vertex-switching
reconstruction problem of Stanley. We study this problem in the language of signed graphs, where the net-degree
of a vertex and the notion of balance come naturally into play. We prove that a signed graph on $n\neq4$ vertices is
reconstructible when its net-degrees are all positive, all negative or all
zero, which generalizes the reconstructibility of regular graphs in Stanley's
problem. For connected balanced signed graphs, whose switching reconstruction problem we
show to be equivalent to that of connected graphs with a distinguished
bipartition of the vertex set, we prove reconstructibility when the two parts of
the bipartition differ in size by more than two, and when all net-degrees are
non-negative.
\end{abstract}

\section{Introduction}

In 2011, Bondy and Mercier \cite[Section~6]{BondyMercier} introduced the switching
reconstruction problem for $2$-edge-colored graphs: the \emph{switching} of such
a graph at a vertex $v$ exchanges the colors of the edges incident with $v$, the
\emph{switching deck} is the multiset of the isomorphism classes of the $n$ graphs
obtained by switching at each vertex, where $n$ is the order of the graph, and the question is whether the graph is
determined, up to isomorphism, by its deck. When the underlying graph is complete
and the two colors are read as ``edge'' and ``non-edge'', this is the
vertex-switching reconstruction problem posed by Stanley in 1985 \cite{Stanley}.

Stanley's problem is a variant of the classical reconstruction problem of Kelly
and Ulam, see \cite{BondyHemminger} for a survey; the latter problem was
considered for balanced signed graphs by Harary \cite{Harary1985}. Stanley asked whether
every graph on at least five vertices is reconstructible; the question is still
open, and appears as a conjecture in Bondy's list of beautiful conjectures in
graph theory \cite{bo}. Stanley proved that graphs on $n\not\equiv0\pmod 4$
vertices are reconstructible, and the conjecture is known for disconnected
\cite{EllinghamRoyle,kra88}, regular and triangle-free graphs \cite{EllinghamRoyle}, for graphs
whose degrees satisfy suitable conditions \cite{kra88,kra96}, and, by computer
search, for $n<16$ \cite{Nie}; moreover, the number of subgraphs with fewer than
$n/2$ vertices isomorphic to a given graph is reconstructible, a switching
analogue of Kelly's Lemma \cite{EllinghamRoyle}. For further results, including
decks obtained by switching several vertices at a time, see
\cite{Ell,kra94,KrasikovRoditty1987,KrasikovRoditty1992,KrasikovRoditty1994}.
Bondy and Mercier \cite{BondyMercier} posed the analogous problems for digraphs,
where non-reconstructible digraphs on eight vertices exist
\cite{BondyMercier,McKaySchweitzer}, and for $2$-edge-colored graphs, observing
that both Stanley's theorem for $n\not\equiv0\pmod 4$ and the analogue of Kelly's
Lemma carry over to the latter. A general axiomatic framework for switching
problems was set up in \cite{Mc19}, for switching-stability, and in the thesis
\cite{Tan}, for reconstruction; among its consequences, a $2$-edge-colored graph
whose underlying graph $\Gamma$ is connected with $|\mathrm{Aut}(\Gamma)|<2^{n/2-2}$
is reconstructible \cite[Theorem~5.18]{Tan}.

In this paper we study the switching reconstruction problem for $2$-edge-colored graphs  in the language of
\emph{signed graphs}: graphs whose edges carry a sign $\pm1$, with switching in
the sense of Zaslavsky \cite{zasgraph}. The two formulations are equivalent, but
they suggest different invariants. The $2$-coloring formulation is symmetric in
the colors; the signed formulation breaks this symmetry, and brings in two
notions: the \emph{net-degree}
of a vertex, i.e.\ the number of positive minus the number of negative edges incident to that vertex, and \emph{balance}, the property that every cycle has an even number of
negative edges (in the social interpretation of signed graphs, the rule ``the
enemy of my enemy is my friend''), with its theory going back to Harary
\cite{Harary} and Zaslavsky \cite{zasgraph,Zaslavsky}. To our knowledge, this
theory has not been used on the switching reconstruction problem before,
although the connection is recorded in Zaslavsky's bibliography
\cite{ZaslavskyBiblio}. Our main results are obtained precisely through these
notions. A signed graph is \emph{reconstructible} if every signed graph with the
same switching deck is isomorphic to it; for a connected balanced signed graph
$(\Gamma,\sigma)$, the \emph{balance partition} is the partition of the vertex set
into two parts such that the negative edges are exactly those joining the two
parts.

\begin{theor}[Theorems~\ref{coro:rego} and~\ref{cor:p1}]\label{thA}
Let $(\Gamma,\sigma)$ be a signed graph on $n\neq4$ vertices. If the net-degrees
of $(\Gamma,\sigma)$ are all positive, all negative, or all zero, then
$(\Gamma,\sigma)$ is reconstructible.
\end{theor}

\begin{theor}[Theorem~\ref{th:nuovo}]\label{thB}
Let $(\Gamma,\sigma)$ be a connected balanced signed graph on $n\neq4$ vertices.
If the two parts of its balance partition differ in cardinality by more than two, then
$(\Gamma,\sigma)$ is reconstructible.
\end{theor}

\begin{theor}[Theorem~\ref{theo:final}]\label{thC}
Let $(\Gamma,\sigma)$ be a connected balanced signed graph on $n\neq4$ vertices.
If every vertex has non-negative net-degree, then $(\Gamma,\sigma)$ is
reconstructible.
\end{theor}

Theorem~\ref{thA} applies in particular to all net-regular signed graphs, those
with constant net-degree, and thus generalizes the reconstructibility of regular
graphs in Stanley's problem \cite{EllinghamRoyle}. The case of net-degrees of one
strict sign rests on a signed version of the \emph{switch partners} of Krasikov
\cite{kra88} and on the identity $d^\pm(v)+d^\pm(w)=2\sigma(v,w)$ relating the
net-degrees of two switch partners (Proposition~\ref{sw.par.form}). Along the
way we prove that the numbers of positive and negative edges and the multiset of
net-degrees of a signed graph on $n\neq4$ vertices are reconstructible
(Proposition~\ref{prop:basic}), and so is a signed graph whose signature is
all-positive or all-negative (Corollary~\ref{cor:all}).

The balanced case is where the signed point of view is most effective. A connected
balanced signed graph is the same thing as a connected graph with a
distinguished bipartition of its vertex set, and switching at $v$ moves $v$ to
the other part (Proposition~\ref{prop:bip}); in the terminology of
\cite[Problem~5]{Mc19}, this is the switching problem for $2$-vertex-colored
graphs up to \emph{weak isomorphism}, i.e.\ isomorphism allowing the global exchange of
the two colors. The sizes of the two parts are reconstructible
(Lemma~\ref{prop:bic}), and Theorem~\ref{thB} says that a non-reconstructible
connected bipartitioned graph must have parts of sizes $\{n/2,n/2\}$ or
$\{n/2-1,n/2+1\}$. Its proof passes through the problem for $2$-vertex-colored
graphs up to color-preserving isomorphism (\cite[Problem~5]{Mc19},
Variant~(a)): when the parts have different sizes they can be told apart in
every card, and an Alon--Coppersmith type argument applies
(Proposition~\ref{teo:2}). For balanced signed graphs, Theorem~\ref{thC} generalizes the non-negative
case of Theorem~\ref{thA}; its core is Lemma~\ref{lem:tec}, which treats the
signed graphs whose net-degrees are $0$ or $2$.  When $\Gamma$ is bipartite, in particular for signed trees,
Theorem~\ref{thC} also applies when all net-degrees are non-positive
(Corollary~\ref{coro:bip}). 

The hypotheses of Theorems~\ref{thA}, \ref{thB} and~\ref{thC} only involve the
net-degrees and the cardinalities of the parts of the balance partition, both of which
are reconstructible, and not the symmetries of $(\Gamma,\sigma)$ or of $\Gamma$; in this respect our
approach is complementary to that of \cite{Tan}. For instance,
Example~\ref{exa:theta} exhibits a signed graph on eight vertices to which
Theorem~\ref{thC} applies, while the results of \cite{Tan} give no information.

The paper is organized as follows. Section~\ref{sec:pre} introduces signed
graphs, switching and the switching deck, and collects what is known from
Stanley's problem. Section~\ref{sec:ns} studies net-degrees and switch partners,
and proves the results that hold for all signed graphs, including
Theorem~\ref{thA}. Section~\ref{sec:balance} is devoted to the balanced case: the
equivalence with bipartitioned graphs, the relation with $2$-vertex-colored
graphs and Theorem~\ref{thB}, and finally Theorem~\ref{thC}. Two natural questions remain open: whether connected balanced signed graphs with
parts of sizes $\{n/2,n/2\}$ or $\{n/2-1,n/2+1\}$ having both negative and
non-negative net-degrees are reconstructible, and whether non-reconstructible
signed graphs on more than four vertices exist at all, a question which
contains Stanley's conjecture when  the underlying graphs are complete.

\section{Preliminaries}\label{sec:pre}

\subsection{Signed graphs and switching}
A \emph{signed graph} is a pair $(\Gamma,\sigma)$, where $\Gamma=(V_\Gamma,E_\Gamma)$ is
a simple graph, the \emph{underlying graph}, and
$\sigma\colon E_\Gamma\to\{-1,+1\}$ is the \emph{signature}. We write $v\sim w$ if
$\{v,w\}\in E_\Gamma$, and $\sigma(v,w)$ for $\sigma(\{v,w\})$; we extend $\sigma$ to
all pairs of distinct vertices by setting $\sigma(v,w)=0$ when $v\nsim w$.
Moreover, $e(\Gamma)=|E_\Gamma|$, and $d(v)$ is the degree of $v$ in $\Gamma$. We
denote by $e^+(\sigma)$ and $e^-(\sigma)$ the numbers of positive and negative
edges, and set $e^\pm(\sigma)=e^+(\sigma)-e^-(\sigma)$. Likewise, $d^+_\sigma(v)$ and
$d^-_\sigma(v)$ are the numbers of positive and negative edges incident to $v$, and
$d^\pm_\sigma(v)=d^+_\sigma(v)-d^-_\sigma(v)$ is the \emph{net-degree} of $v$. We write
$\deg(\Gamma)$ for the multiset of the degrees $d(v)$ and $\deg^\pm(\Gamma,\sigma)$ for
the multiset of the net-degrees $d^\pm_\sigma(v)$, $v\in V_\Gamma$. A \emph{signed
subgraph} of $(\Gamma,\sigma)$ is a subgraph of $\Gamma$ with the restriction of $\sigma$.

Two signed graphs $(\Gamma_1,\sigma_1)$ and $(\Gamma_2,\sigma_2)$ are
\emph{isomorphic}, written $(\Gamma_1,\sigma_1)\cong(\Gamma_2,\sigma_2)$, if there is a
bijection $\phi\colon V_{\Gamma_1}\to V_{\Gamma_2}$ with
$\sigma_1(u,v)=\sigma_2(\phi(u),\phi(v))$ for all distinct $u,v\in V_{\Gamma_1}$. The
\emph{switching} of $\sigma$ at a vertex $v$ is the signature $\sigma_v$ obtained by
changing the sign of every edge incident to $v$ \cite{zasgraph}. Switchings commute, and
for $A\subseteq V_\Gamma$ we write $\sigma_A$ for the signature obtained by switching at
every vertex of $A$ in any order, i.e.\ by changing the sign of every edge with
exactly one endpoint in $A$; for $A=\{v,w\}$ we simply write $\sigma_{vw}$. Two
signatures $\sigma,\sigma'$ on $\Gamma$ are \emph{switching equivalent} if
$\sigma'=\sigma_A$ for some $A$, and two signed graphs are \emph{switching isomorphic}
if one is isomorphic to a switching of the other.

\begin{remark}\label{rem:sei}
A graph $\Gamma$ on $n$ vertices determines the signed complete graph
$(K_n,\sigma_\Gamma)$ on the same vertex set whose negative edges are the edges of
$\Gamma$ (one may equally use the positive edges, replacing $\Gamma$ by its
complement). Under this correspondence, Seidel switching of $\Gamma$ at $v$, which
exchanges the neighbours and the non-neighbours of $v$, is the switching of
$(K_n,\sigma_\Gamma)$ at $v$ \cite{sei,zasgraph}.
\end{remark}

A signature $\sigma$, and the signed graph $(\Gamma,\sigma)$, is \emph{balanced} if
every cycle of $\Gamma$ contains an even number of negative edges. Equivalently,
$V_\Gamma$ can be partitioned into two sets $A,B$ such that the negative edges are
exactly those joining $A$ and $B$ \cite{Harary}, or $\sigma$ is switching
equivalent to the all-positive signature \cite{zasgraph}.

\subsection{The switching deck}
The \emph{switching deck} of a signed graph $(\Gamma,\sigma)$ on $n$ vertices is the
multiset of the $n$ isomorphism classes of the signed graphs $(\Gamma,\sigma_v)$,
$v\in V_\Gamma$, which are called the \emph{cards}. Isomorphic signed graphs have the
same deck. An isomorphism invariant of signed graphs is \emph{reconstructible} if
it is determined by the deck, i.e.\ if it takes the same value on any two signed
graphs with the same deck. The signed graph $(\Gamma,\sigma)$ is
\emph{reconstructible} if every signed graph with the same deck is isomorphic to
it, that is, if its isomorphism class is reconstructible.
\begin{problem}\label{pro:sig}
Which signed graphs are determined, up to isomorphism, by their switching deck?
More generally, which properties of a signed graph are reconstructible?
\end{problem}
Since changing all signs commutes with switching and with isomorphisms, $(\Gamma,\sigma)$ is reconstructible if and only if $(\Gamma,-\sigma)$ is.
\begin{figure}
\begin{center}
\begin{tikzpicture}[scale=0.75,
    edgepos/.style={solid},
    edgeneg/.style={dashed},
    vertex/.style={circle,fill,inner sep=1pt}
]
\draw (0,0) rectangle (4,-8);
\foreach \y in {-2,-4,-6} \draw (0,\y) -- (4,\y);

\def\pawcoords{
  \node[vertex] (a) at (0,1) {};
  \node[vertex] (b) at (0,0) {};
  \node[vertex] (c) at (-1,-1) {};
  \node[vertex] (d) at (1,-1) {};
}

\node at (0.3,-0.3) {1};
\begin{scope}[shift={(2,-1)},scale=0.7]
  \pawcoords
  \draw[edgeneg] (a)--(b);
  \draw[edgeneg] (b)--(c);
  \draw[edgepos] (b)--(d);
  \draw[edgeneg] (c)--(d);
\end{scope}

\node at (0.3,-2.3) {2};
\begin{scope}[shift={(2,-3)},scale=0.7]
  \pawcoords
  \draw[edgeneg] (a)--(b);
  \draw[edgepos] (b)--(c);
  \draw[edgeneg] (b)--(d);
  \draw[edgeneg] (c)--(d);
\end{scope}

\node at (0.3,-4.3) {3};
\begin{scope}[shift={(2,-5)},scale=0.7]
  \pawcoords
  \draw[edgepos] (a)--(b);
  \draw[edgepos] (b)--(c);
  \draw[edgepos] (b)--(d);
  \draw[edgepos] (c)--(d);
\end{scope}

\node at (0.3,-6.3) {4};
\begin{scope}[shift={(2,-7)},scale=0.7]
  \pawcoords
  \draw[edgepos] (a)--(b);
  \draw[edgeneg] (b)--(c);
  \draw[edgeneg] (b)--(d);
  \draw[edgepos] (c)--(d);
\end{scope}

\begin{scope}[shift={(-4,-4)}]
  \pawcoords
  \node at (0,1.3) {\tiny 1};
  \node at (0.3,0) {\tiny 2};
  \node at (-1.3,-1) {\tiny 3};
  \node at (1.3,-1) {\tiny 4};
  \draw[edgepos] (a)--(b);
  \draw[edgeneg] (b)--(c);
  \draw[edgepos] (b)--(d);
  \draw[edgeneg] (c)--(d);
\end{scope}

\begin{scope}[shift={(8,-4)}]
  \pawcoords
  \node at (0,1.3) {\tiny 4};
  \node at (0.3,0) {\tiny 3};
  \node at (-1.3,-1) {\tiny 2};
  \node at (1.3,-1) {\tiny 1};
  \draw[edgeneg] (a)--(b);
  \draw[edgeneg] (b)--(c);
  \draw[edgeneg] (b)--(d);
  \draw[edgepos] (c)--(d);
\end{scope}

\end{tikzpicture}
\end{center}
\caption{Two signed graphs (on the left and on the right)  with the same switching deck (in the center).}\label{fig:nonrec}
\end{figure}
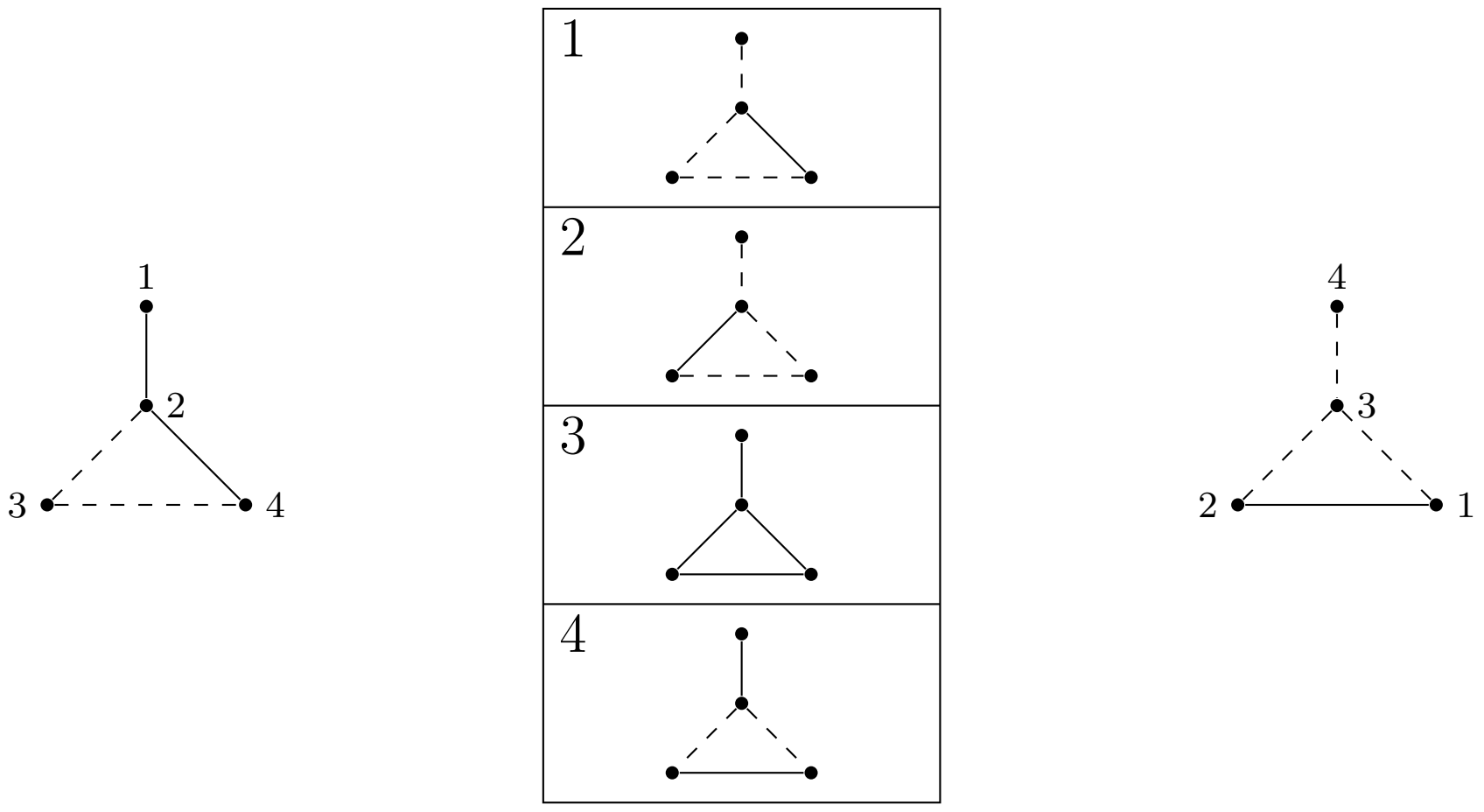

\begin{example}\label{ex:ric}
Fig.~\ref{fig:nonrec} shows two non-isomorphic signed graphs on four vertices with
the same switching deck. Here and in all figures, solid edges are positive and
dashed edges negative. One checks directly that exactly five isomorphism
classes of signed graphs on four vertices are reconstructible, as noted in
\cite{BondyMercier}.
\end{example}

Whether non-reconstructible signed graphs exist on more than four vertices is
open. In \cite{BondyMercier} the question was raised for eight vertices; in
\cite{Tan} the existence of such examples on eight vertices is claimed, but no
construction or proof is provided.

By Remark~\ref{rem:sei}, the vertex-switching reconstruction problem of Stanley
\cite{Stanley} is the case of Problem~\ref{pro:sig} in which the underlying graph
is complete, and every result on Stanley's problem applies to signed complete
graphs. Thus $(K_n,\sigma)$, $n\neq4$, is reconstructible when its negative (or
positive) edges form a disconnected \cite{kra88,EllinghamRoyle}, regular or
triangle-free \cite{EllinghamRoyle} graph, when their degrees lie outside certain
ranges \cite{kra88,kra96}, and for $n<16$ \cite{Nie}. One result, however, extends
to all signed graphs: Stanley's theorem that graphs on $n\not\equiv0 \pmod 4$
vertices are reconstructible, whose  proof (see
\cite[Theorem~1.1]{BondyMercier}) applies verbatim to $2$-edge-colored graphs,
hence to signed graphs, as observed in \cite[Section~6]{BondyMercier} and, in a
general framework, in \cite{Tan}.

\begin{proposition}\label{prop:4}
Let $(\Gamma,\sigma)$ be a signed graph on $n$ vertices. If $n\not\equiv0 \pmod 4$,
then $(\Gamma,\sigma)$ is reconstructible.
\end{proposition}

\section{Net-degrees and switch partners}\label{sec:ns}
\subsection{Basic reconstructible invariants}\label{sec:basic}
Every card of the switching deck of $(\Gamma,\sigma)$ determines the isomorphism
class of $\Gamma$ and the switching isomorphism class of $(\Gamma,\sigma)$. Hence
all invariants of $\Gamma$, such as the order $n$, the number of edges $e(\Gamma)$,
the degree multiset $\deg(\Gamma)$, connectedness and the number of triangles,
and all switching invariants of $(\Gamma,\sigma)$, such as balance and spectrum,
are reconstructible.

\begin{proposition}\label{prop:basic}
Let $(\Gamma,\sigma)$ be a signed graph on $n\neq 4$ vertices.
\begin{enumerate}[(i)]
\item $e^\pm(\sigma)$, $e^+(\sigma)$ and $e^-(\sigma)$ are reconstructible.
\item Each card determines the net-degree of the vertex at which it was obtained.
\item The multiset $\deg^\pm(\Gamma,\sigma)$ of net-degrees is reconstructible.
\end{enumerate}
\end{proposition}
\begin{proof}
Switching at $v$ turns $d^+_\sigma(v)$ positive edges into negative ones and vice
versa, so
\begin{equation}\label{eq:epm}
e^\pm(\sigma_v)=e^\pm(\sigma)-2d^\pm_\sigma(v).
\end{equation}
Summing over $v$ and using $\sum_v d^\pm_\sigma(v)=2e^\pm(\sigma)$ gives
$\sum_v e^\pm(\sigma_v)=(n-4)\,e^\pm(\sigma)$. The left side can be computed from
the deck, so $e^\pm(\sigma)$ is reconstructible,
and so are $e^+(\sigma)=\tfrac12\bigl(e(\Gamma)+e^\pm(\sigma)\bigr)$ and
$e^-(\sigma)=\tfrac12\bigl(e(\Gamma)-e^\pm(\sigma)\bigr)$, which proves (i).
By (i) and
Eq.~\eqref{eq:epm}, $d^\pm_\sigma(v)=\tfrac12\bigl(e^\pm(\sigma)-e^\pm(\sigma_v)\bigr)$
is determined by the card $(\Gamma,\sigma_v)$, which proves (ii), and (iii)
follows.
\end{proof}

\begin{corollary}\label{cor:all}
If $n\neq 4$ and $\sigma$ is all-positive or all-negative, then $(\Gamma,\sigma)$ is
reconstructible.
\end{corollary}

The switching analogue of Kelly's Lemma holds  for signed graphs.
\begin{proposition}[\cite{EllinghamRoyle}; see {\cite[Cor.~5.1, Sec.~6]{BondyMercier}}]\label{prop:kelly}
If $(F,\tau)$ has fewer than $n/2$ vertices, the number of signed subgraphs of
$(\Gamma,\sigma)$ isomorphic to $(F,\tau)$ is reconstructible.
\end{proposition}
The conditions of Proposition~\ref{prop:kelly} are not sufficient for two signed
graphs to have the same deck: see Example~\ref{exa:due}  in
Section~\ref{sec:balance}.

\subsection{Net-regular signed graphs}\label{sec:netreg}
A signed graph $(\Gamma,\sigma)$ is \emph{net-regular} if $d^\pm_\sigma(v)=d$ for all
$v\in V_\Gamma$ and some $d\in\mathbb Z$. A complete signed graph is net-regular if
and only if its positive (equivalently, negative) edges form a regular graph, so
that Theorem~\ref{coro:rego} below generalizes \cite[Theorem~4.5]{EllinghamRoyle}
via the correspondence of Remark~\ref{rem:sei}. We first need a simple lemma.

\begin{lemma}\label{lemma:intorno}
Let $u,v$ be vertices of $(\Gamma,\sigma)$ such that $\sigma(u,w)=\sigma(v,w)$ for
every $w\notin\{u,v\}$. Then $(\Gamma,\sigma_u)\cong(\Gamma,\sigma_v)$.
\end{lemma}
\begin{proof}
Under this hypothesis, the transposition $\tau$ of $u$ and $v$ is an isomorphism
from $(\Gamma,\sigma_u)$ to $(\Gamma,\sigma_v)$, since
$\sigma_v(\tau(u),\tau(w))=-\sigma(v,w)=-\sigma(u,w)=\sigma_u(u,w)$ for every
$w\notin\{u,v\}$, and the signs of the other pairs are not affected.
\end{proof}

\begin{theorem}\label{coro:rego}
Let $n\neq4$. If $(\Gamma,\sigma)$ is net-regular, then it is reconstructible.
\end{theorem}
\begin{proof}
Let $d$ be the common net-degree, which is reconstructible by
Proposition~\ref{prop:basic}(iii), and fix a card $(\Gamma,\sigma_v)$. In $\sigma_v$
the vertex $v$ has net-degree $-d$, a vertex adjacent to $v$ by a positive
(resp.\ negative) edge of $\sigma_v$ has net-degree $d+2$ (resp.\ $d-2$), and every
other vertex has net-degree $d$. Set
$$
I^+=\{x\in V_\Gamma: d^\pm_{\sigma_v}(x)=d+2\},\qquad
I^-=\{x\in V_\Gamma: d^\pm_{\sigma_v}(x)=d-2\}.
$$
Both sets are determined by the card, and the positive (resp.\ negative)
neighbours of $v$ in $\sigma_v$ are exactly the vertices of $I^+\setminus\{v\}$
(resp.\ $I^-\setminus\{v\}$). Now choose any vertex $v'$ whose positive and
negative neighbours in $\sigma_v$ are exactly the vertices of $I^+\setminus\{v'\}$
and of $I^-\setminus\{v'\}$; such a $v'$ exists, since $v$ has this property, and
it is determined by the card. Then $v$ and $v'$ satisfy the hypothesis of
Lemma~\ref{lemma:intorno} in $(\Gamma,\sigma_v)$, hence
$(\Gamma,\sigma_{vv'})\cong(\Gamma,\sigma_{vv})=(\Gamma,\sigma)$, and $(\Gamma,\sigma)$
is reconstructible.
\end{proof}
Note that a single card suffices, provided net-regularity is known a priori; and
if $d\notin\{-1,0,1\}$ then $v$ is simply the unique vertex of net-degree $-d$ in
the card.
\begin{remark}
Theorem~\ref{coro:rego} refines the observation of
\cite[Section~6]{BondyMercier} that a graph or digraph whose switchings are all
isomorphic is switching-reconstructible: for $n\neq4$ it is enough that all cards
have the same number of positive edges, since then, by
Proposition~\ref{prop:basic}(ii), all vertices have the same net-degree.
\end{remark}

\subsection{Switch partners}\label{sec:sp}
The following notion goes back to Krasikov~\cite[Lemma~1]{kra88} in the setting of
Stanley's problem; the name is due to Niesink~\cite{Nie}.

\begin{definition}
Let $v$ be a vertex of a signed graph $(\Gamma,\sigma)$. A vertex $w\neq v$ is a
\emph{switch partner} of $v$ if $(\Gamma,\sigma_{vw})$ is not isomorphic to
$(\Gamma,\sigma)$ and has the same switching deck.
\end{definition}

\begin{example}\label{ex:sp}
In the signed graph on the right of Fig.~\ref{fig:nonrec}, both $1$ and $2$ are
switch partners of $4$: switching at $4$ and then at $1$, or at $4$ and then at $2$,
yields a signed graph isomorphic to the one on the left, which has the same deck.
In particular switch partners need not be unique.
\end{example}

\begin{proposition}\label{lem:swp}
Let $(\Gamma,\sigma)$ be non-reconstructible, and let $\sigma'$ be a signature on
$\Gamma$ with the same switching deck and $(\Gamma,\sigma')\not\cong(\Gamma,\sigma)$.
Then every $v\in V_\Gamma$ has a switch partner $w$ with
$(\Gamma,\sigma_{vw})\cong(\Gamma,\sigma')$.
\end{proposition}
\begin{proof}
Since the decks coincide, for every $v$ there are $z\in V_\Gamma$ and an isomorphism
$\phi$ from $(\Gamma,\sigma_v)$ to $(\Gamma,\sigma'_z)$. An isomorphism of signed graphs
commutes with switching, i.e.\ $\phi$ is also an isomorphism
$(\Gamma,\sigma_{vu})\to(\Gamma,\sigma'_{z\phi(u)})$ for every $u$. Taking
$w=\phi^{-1}(z)$ we get $(\Gamma,\sigma_{vw})\cong(\Gamma,\sigma'_{zz})=(\Gamma,\sigma')$,
and $w\neq v$ because $\sigma'\not\cong\sigma=\sigma_{vv}$.
\end{proof}

The next identity is the signed analogue of \cite[Corollary~1]{kra88} and is a crucial
tool for the rest of the paper.

\begin{proposition}\label{sw.par.form}
Let $(\Gamma,\sigma)$ be a signed graph on $n\neq4$ vertices and let $w$ be a switch
partner of $v$. Then
$$
d^\pm_\sigma(v)+d^\pm_\sigma(w)=2\,\sigma(v,w).
$$
\end{proposition}
\begin{proof}
By Proposition~\ref{prop:basic}(i), $e^\pm(\sigma_{vw})=e^\pm(\sigma)$, since the two
signed graphs have the same deck. Applying Eq.~\eqref{eq:epm} first to $\sigma$ at $v$
and then to $\sigma_v$ at $w$ gives
$e^\pm(\sigma_{vw})=e^\pm(\sigma)-2d^\pm_\sigma(v)-2d^\pm_{\sigma_v}(w)$, hence
$d^\pm_\sigma(v)+d^\pm_{\sigma_v}(w)=0$. Finally
$d^\pm_{\sigma_v}(w)=d^\pm_\sigma(w)-2\sigma(v,w)$, because only the edge
$\{v,w\}$, if present, changes sign at $w$ when switching at $v$.
\end{proof}
Proposition~\ref{sw.par.form} leaves only three possible net-degrees for a switch
partner, which gives a first criterion.
\begin{proposition}\label{coro:p}
Let $n\neq 4$. If some $d\in\deg^\pm(\Gamma,\sigma)$ satisfies
$\{-d-2,-d,-d+2\}\cap\deg^\pm(\Gamma,\sigma)=\emptyset$, then $(\Gamma,\sigma)$ is
reconstructible.
\end{proposition}
\begin{proof}
By Proposition~\ref{sw.par.form}, a vertex of net-degree $d$ has no switch
partner, so $(\Gamma,\sigma)$ is reconstructible by Proposition~\ref{lem:swp}.
\end{proof}
The condition is easiest to check on the largest net-degree.
\begin{proposition}\label{cor:d1dn}
Let $n\neq4$ and let $d_1>\dots>d_k$ be the distinct net-degrees of $(\Gamma,\sigma)$.
If $d_1>-d_k+2$, then $(\Gamma,\sigma)$ is reconstructible.
\end{proposition}
\begin{proof}
The hypothesis gives $\{-d_1-2,-d_1,-d_1+2\}\cap\deg^\pm(\Gamma,\sigma)=\emptyset$,
and Proposition~\ref{coro:p} applies with $d=d_1$.
\end{proof}

\begin{remark}
For $\Gamma=K_n$ with positive edges forming a graph of minimum degree $\delta$ and
maximum degree $\Delta$ (Remark~\ref{rem:sei}), one has $d_1=2\Delta-n+1$ and
$d_k=2\delta-n+1$, so Proposition~\ref{cor:d1dn} says that a non-reconstructible
graph satisfies $\delta+\Delta\le n$ and, applied to the complementary graph,
$\delta+\Delta\ge n-2$: this is \cite[Corollary~2]{kra88}.
\end{remark}

\begin{theorem}\label{cor:p1}
Let $n\neq4$. If all net-degrees of $(\Gamma,\sigma)$ are positive, or all are
negative, then $(\Gamma,\sigma)$ is reconstructible.
\end{theorem}
\begin{proof}
Replacing $\sigma$ by $-\sigma$ we may assume $d_1>\dots>d_k>0$. If $d_1=1$ the
signed graph is net-regular and Theorem~\ref{coro:rego} applies. If $d_1\ge2$ then
$-d_1+2\le0$, so Proposition~\ref{coro:p} applies with $d=d_1$.
\end{proof}
\begin{remark}
For $\Gamma=K_n$, by Remark~\ref{rem:sei}, Theorem~\ref{cor:p1} says that a graph
with minimum degree at least $n/2$, or maximum degree less than $n/2$, is
reconstructible in Stanley's problem. This also follows from
\cite[Corollary~2]{kra88} together with the reconstructibility of regular graphs
\cite[Theorem~4.5]{EllinghamRoyle}.
\end{remark}

\section{Switching reconstruction of balanced signed graphs}\label{sec:balance}
Recall that $(\Gamma,\sigma)$ is balanced if and only if $V_\Gamma$ admits a partition
$\{A,B\}$ such that the negative edges are exactly those joining $A$ and $B$
\cite{Harary}. When $\Gamma$ is connected this partition is unique, so a connected
balanced signed graph is the same thing as a connected graph with a distinguished
bipartition of its vertex set. We first make this precise, then use the bipartition to obtain
reconstructibility results specific to the balanced case.

\subsection{Bipartitioned graphs}\label{sec:bip}
A \emph{bipartitioned graph} is a pair $(\Gamma,\pi)$ where $\pi=\{A,B\}$ is a
partition of $V_\Gamma$ into two (possibly empty) parts; we write $v\sim_\pi w$ if
$v,w$ lie in the same part. Two bipartitioned graphs $(\Gamma,\pi)$, $(\Gamma',\pi')$
are \emph{isomorphic} if there is a graph isomorphism $\phi\colon\Gamma\to\Gamma'$
with $v\sim_\pi w \iff \phi(v)\sim_{\pi'}\phi(w)$; note that $\phi$ may interchange
the two parts. The \emph{switching} of $(\Gamma,\pi)$ at $v$ is $(\Gamma,\pi_v)$,
where $\pi_v$ is obtained from $\pi$ by moving $v$ to the other part; for
$W\subseteq V_\Gamma$, $\pi_W$ is defined accordingly. The \emph{switching deck} of
$(\Gamma,\pi)$ is the multiset of the isomorphism classes of the $n$ switchings
$(\Gamma,\pi_v)$, and reconstructibility is defined as in Section~\ref{sec:pre}.

\begin{problem}\label{pro:bip}
Which bipartitioned graphs are determined, up to isomorphism, by their switching
deck?
\end{problem}
Bipartitioned graphs, with the switching defined above, are exactly the
$2$-vertex-colored graphs up to \emph{weak isomorphism}, i.e.\ isomorphism
allowing the global exchange of the two colors, of \cite[Problem~5,
Variant~(b)]{Mc19}, where those whose switchings are all isomorphic to the
original are classified.

Given $(\Gamma,\pi)$, let $\sigma_\pi$ be the signature of $\Gamma$ that is positive
on edges inside a part and negative on edges between the parts.

\begin{proposition}\label{prop:bip}
Let $\Gamma$ be a connected graph. The map $\pi\mapsto\sigma_\pi$ is a bijection
between the bipartitions of $V_\Gamma$ and the balanced signatures of $\Gamma$, and
$\sigma_{\pi_v}=(\sigma_\pi)_v$ for every $v\in V_\Gamma$. Moreover, if $\Gamma'$ is
connected and $\pi'$ is a bipartition of $V_{\Gamma'}$, then
$(\Gamma,\pi)\cong(\Gamma',\pi')$ if and only if
$(\Gamma,\sigma_\pi)\cong(\Gamma',\sigma_{\pi'})$. Consequently, $(\Gamma,\pi)$ is
reconstructible if and only if $(\Gamma,\sigma_\pi)$ is.
\end{proposition}
\begin{proof}
Every $\sigma_\pi$ is balanced, and every balanced signature arises in this way
\cite{Harary}. Since $\Gamma$ is connected, $v\sim_\pi w$ holds if and only if some
(equivalently, every) path from $v$ to $w$ carries an even number of negative edges
of $\sigma_\pi$; so $\pi$ is determined by $\sigma_\pi$, and the map is a bijection.
Moving $v$ to the other part changes exactly the signs of the edges incident to $v$, so
$\sigma_{\pi_v}=(\sigma_\pi)_v$. The same path argument shows that a graph
isomorphism $\phi\colon\Gamma\to\Gamma'$ satisfies
$\sigma_\pi(v,w)=\sigma_{\pi'}(\phi(v),\phi(w))$ for all adjacent $v,w$ if and only if
$v\sim_\pi w\iff\phi(v)\sim_{\pi'}\phi(w)$ for all $v,w$, that is, $\phi$ is an
isomorphism of signed graphs if and only if it is an isomorphism of bipartitioned
graphs. Hence the deck of $(\Gamma,\pi)$ corresponds to the deck of
$(\Gamma,\sigma_\pi)$, and since underlying graph and balance are reconstructible,
$(\Gamma,\pi)$ is reconstructible if and only if $(\Gamma,\sigma_\pi)$ is.
\end{proof}
For connected graphs, Problem~\ref{pro:bip} is therefore the restriction of
Problem~\ref{pro:sig} to balanced signed graphs, and the results of
Section~\ref{sec:ns} translate immediately, as in the following corollary. 
\begin{corollary}\label{cor:bipbasic}
Let $(\Gamma,\pi)$ be a connected bipartitioned graph on $n\neq4$ vertices. Then
$(\Gamma,\pi)$ is reconstructible in each of the following cases:
\begin{enumerate}[(i)]
\item $n\not\equiv0 \pmod 4$;
\item $\pi$ is the trivial partition $\{V_\Gamma,\emptyset\}$;
\item every vertex has as many neighbours in its own part as in the other part.
\end{enumerate}
\end{corollary}
\begin{proof}
By Proposition~\ref{prop:bip}, apply Proposition~\ref{prop:4},
Corollary~\ref{cor:all} (the signature $\sigma_\pi$ is all-positive) and
Theorem~\ref{coro:rego} ($(\Gamma,\sigma_\pi)$ is net-regular of net-degree $0$),
respectively.
\end{proof}
For a connected balanced $(\Gamma,\sigma)$ we denote by $\pi_\sigma$ the unique
partition with $\sigma_{\pi_\sigma}=\sigma$, the \emph{balance partition}, and by
$\pi_\sigma(v)$ the part containing $v$. Complementation behaves well with respect
to this partition. In \cite{CDSdm} a balanced signature $\sigma^c$ on $\Gamma^c$ is
associated with any signed graph $(\Gamma,\sigma)$; when $(\Gamma,\sigma)$ is
connected and balanced, $\sigma^c$ is simply $\sigma_{\pi_\sigma}$ computed on
$\Gamma^c$, i.e.\ the balanced signature of $\Gamma^c$ with the same balance
partition.

\begin{proposition}\label{prop:comple}
Let $(\Gamma,\sigma)$ be a connected balanced signed graph with $\Gamma^c$
connected. Then $(\Gamma,\sigma)$ is reconstructible if and only if
$(\Gamma^c,\sigma^c)$ is.
\end{proposition}
\begin{proof}
The cards of $(\Gamma^c,\pi)$ are the complements of the cards of $(\Gamma,\pi)$,
so $(\Gamma,\pi)$ is reconstructible if and only if $(\Gamma^c,\pi)$ is; apply
Proposition~\ref{prop:bip} on both sides.
\end{proof}

\begin{lemma}\label{prop:bic}
Let $(\Gamma,\pi)$ be a connected bipartitioned graph on $n\neq4$ vertices,
$\pi=\{A,B\}$. The multiset $\{|A|,|B|\}$ is reconstructible.
\end{lemma}
\begin{proof}
Let $(\Gamma,\pi')$, $\pi'=\{A',B'\}$, have the same deck as $(\Gamma,\pi)$; we show
that $\{|A'|,|B'|\}=\{|A|,|B|\}$. If one of the two partitions is trivial, the
corresponding bipartitioned graph is reconstructible by
Corollary~\ref{cor:bipbasic}(ii), so the two are isomorphic and the claim holds.
Let then $a=|A|$, $b=|B|$ with $a,b\ge1$, and likewise for $\pi'$. The deck
contains $a$ cards whose parts have sizes $\{a-1,b+1\}$ and $b$ cards whose parts
have sizes $\{a+1,b-1\}$. If $a\neq b$ these two multisets are different, and
counting the cards of each kind gives $\{a,b\}$; if $a=b$, all cards have parts
of sizes $\{n/2-1,n/2+1\}$, which cannot happen when $a\neq b$. Hence
$\{a,b\}$ is determined by the deck, and the claim follows.
\end{proof}

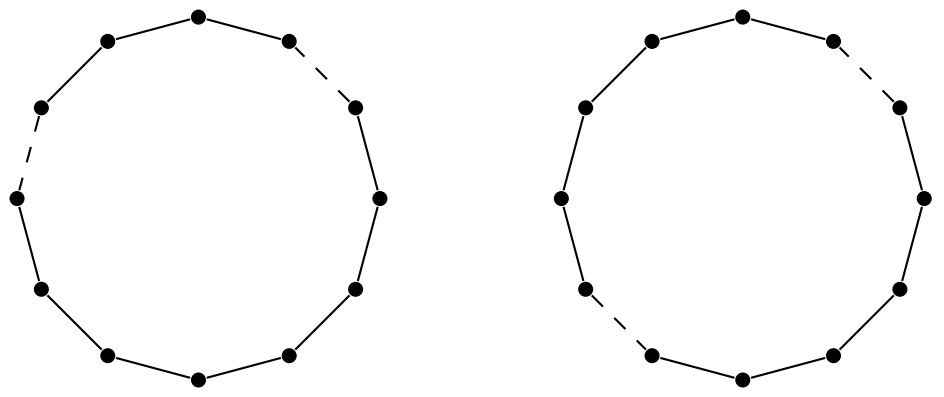
\begin{figure}
\centering

\begin{tikzpicture}[scale=0.6]

\def\n{12}
\def\r{2cm}

\foreach \i in {1,...,\n} {
  \node[fill=black, circle, inner sep=1pt] (A\i) at ({360/\n*(\i - 1) + 90}:\r) {};
}

\foreach \i in {1,...,12} {
  \pgfmathtruncatemacro{\j}{mod(\i,12)+1}
  \ifnum\i=3
    \draw[dashed] (A\i) -- (A\j); 
  \else\ifnum\i=11
    \draw[dashed] (A\i) -- (A\j); 
  \else
    \draw (A\i) -- (A\j);
  \fi\fi
}
\node at (0,-2.8) {$(C_{12}, \sigma_1)$};

\begin{scope}[xshift=6cm]
\foreach \i in {1,...,\n} {
  \node[fill=black, circle, inner sep=1pt] (B\i) at ({360/\n*(\i - 1) + 90}:\r) {};
}

\foreach \i in {1,...,12} {
  \pgfmathtruncatemacro{\j}{mod(\i,12)+1}
  \ifnum\i=5
    \draw[dashed] (B\i) -- (B\j); 
  \else\ifnum\i=11
    \draw[dashed] (B\i) -- (B\j); 
  \else
    \draw (B\i) -- (B\j);
  \fi\fi
}
\node at (0,-2.8) {$(C_{12}, \sigma_2)$};
\end{scope}
\end{tikzpicture}
\caption{
The signed graphs $(C_{12}, \sigma_1)$ and $(C_{12}, \sigma_2)$ in Example~\ref{exa:due}.}\label{fig:cicli}
\end{figure}
\begin{example}\label{exa:due}
Let $v_1,\dots,v_{12}$ be the vertices of $C_{12}$ in cyclic order, and let
$\sigma_1,\sigma_2$ be the balanced signatures whose negative edges are
$\{v_3,v_4\},\{v_{11},v_{12}\}$ and $\{v_5,v_6\},\{v_{11},v_{12}\}$ respectively,
shown in Fig.~\ref{fig:cicli}.
A direct check shows that, for every signed graph $(F,\tau)$ on at most $5$
vertices, $(C_{12},\sigma_1)$ and $(C_{12},\sigma_2)$ contain the same number of
signed subgraphs isomorphic to $(F,\tau)$, so Proposition~\ref{prop:kelly} does not
separate their decks. Their balance partitions have sizes $\{4,8\}$ and $\{6,6\}$,
so Lemma~\ref{prop:bic} does.
\end{example}

\begin{remark}
Lemma~\ref{prop:bic} is equivalent to the reconstructibility of balanced complete
signed graphs on $n\neq4$ vertices, via the balanced completion of \cite{CDSdm};
the latter is known, since by Remark~\ref{rem:sei} it amounts to the
reconstructibility of disconnected graphs in Stanley's problem \cite{kra88}.
\end{remark}

We now turn to a closely related problem, which will provide the main tool of this
subsection. A \emph{$2$-vertex-colored graph} is a pair $(\Gamma,c)$ with
$c\colon V_\Gamma\to\{0,1\}$ an arbitrary (not necessarily proper) coloring;
isomorphisms are color-preserving graph isomorphisms, the switching at $v$ changes
the color of $v$ only, and the switching deck and reconstructibility are defined as
above. This problem is \cite[Problem~5, Variant~(a)]{Mc19}, where its
switching-stability version is solved. A bipartitioned graph is a $2$-vertex-colored graph in which isomorphisms are
also allowed to exchange the two colors globally. Neither problem reduces to
the other: the two monochromatic colorings of $K_2$ have the same deck in both
senses, and are isomorphic as bipartitioned graphs but not as colored graphs;
the monochromatic $C_4$ and a $C_4$ with two vertices of each color are not
isomorphic in either sense, and have the same bipartitioned deck but different
colored decks.

\begin{lemma}\label{lem:01}
Let $(\Gamma,c)$ be a $2$-vertex-colored graph on $n\neq2$ vertices with $c_0$
black and $c_1$ white vertices. Then $c_0$ and $c_1$ are reconstructible.
\end{lemma}
\begin{proof}
If $c_0\in\{0,n\}$ all cards have $1$ or $n-1$ black vertices, and the two cases
are distinguished since $n\neq2$. Otherwise the deck has $c_0$ cards with $c_0-1$
black vertices and $n-c_0$ cards with $c_0+1$ black vertices, so $c_0$ is the number of cards with
the smaller number of black vertices.
\end{proof}

The next statement can be derived from the general framework of
\cite[Theorem~5.16]{Tan}, combined with our Lemma~\ref{lem:01}. Since our use of it
is limited to the case $c_0\neq c_1$, we give a direct and self-contained proof,
which refines the Alon--Coppersmith argument (see
\cite[Theorem~1.1]{BondyMercier}) by means of a filtering of the deck.
\begin{proposition}\label{teo:2}
Let $(\Gamma,c)$ be a $2$-vertex-colored graph on $n>2$ vertices with
$c_0\neq c_1$. Then $(\Gamma,c)$ is reconstructible.
\end{proposition}
\begin{proof}
We represent multisets of isomorphism classes of $2$-vertex-colored graphs as
formal sums with non-negative integer coefficients, and we denote by
$\mathcal{C}$ the isomorphism class of $(\Gamma,c)$. For $k\in\{0,1,\ldots,n\}$,
let $D_k(\mathcal{C})$ be the formal sum of the $\binom nk$ isomorphism classes
obtained from $\mathcal{C}$ by switching $k$ distinct vertices; thus
$D_0(\mathcal{C})=\mathcal{C}$ and $D_1(\mathcal{C})$ is the switching deck of
$(\Gamma,c)$. The operation $D_k$ extends by linearity to formal sums. Counting
the sequences of $k+1$ switchings in two ways one gets
$D_k(D_1(\mathcal{C}))=D_1(D_k(\mathcal{C}))=(n-k+1)D_{k-1}(\mathcal{C})+(k+1)D_{k+1}(\mathcal{C})$
(see the proof of \cite[Theorem~1.1]{BondyMercier}), and since
$D_k(D_1(\mathcal{C}))$ is computed from the deck,
\begin{equation}\label{eq:ac}
(n-k+1)D_{k-1}(\mathcal{C})+(k+1)D_{k+1}(\mathcal{C}) \mbox{ is reconstructible for every $k=1,\ldots,n-1$.}
\end{equation}

If $n$ is odd, we conclude as in the Alon--Coppersmith argument: by
Eq.~\eqref{eq:ac} (with $k=1$), $n\mathcal{C}+2D_2(\mathcal{C})$ is reconstructible,
and $\mathcal{C}$ is the only class appearing in it an odd number of times, hence
$\mathcal{C}$ is reconstructible. We may therefore assume that $n$ is even.

Let $c_0$ and $c_1$ be the numbers of black and white vertices of $\mathcal{C}$,
which are reconstructible by Lemma~\ref{lem:01}. For a formal sum $D$ of
isomorphism classes, let $\overline{D}$ be the sub-sum of the classes having
exactly $c_0$ black and $c_1$ white vertices. Clearly
$\overline{D_0(\mathcal{C})}=\mathcal{C}$; more generally,
$\overline{D_k(\mathcal{C})}$ is the sum over those sets of $k$ switched vertices
that contain as many black as white vertices. The operation is linear, and if
$D$ is reconstructible then so is $\overline{D}$; combining with
Eq.~\eqref{eq:ac},
\begin{equation}\label{eq:acover}
(n-k+1)\overline{D_{k-1}(\mathcal{C})}+(k+1)\overline{D_{k+1}(\mathcal{C})} \mbox{ is reconstructible for every $k=1,\ldots,n-1$.}
\end{equation}

Exchanging the colors, we may assume $c_0<c_1$. Since $n=c_0+c_1$ is even,
$c_0$ and $c_1$ are not consecutive, hence $2c_0+2\le n$. By
Eq.~\eqref{eq:acover} with $k=2c_0+1$, one has that
$(n-2c_0)\,\overline{D_{2c_0}(\mathcal{C})}+(2c_0+2)\,\overline{D_{2c_0+2}(\mathcal{C})}$
is reconstructible. But $\overline{D_{2c_0+2}(\mathcal{C})}$ is empty, since a set
of $2c_0+2$ switched vertices selected by the overline would contain $c_0+1$
black vertices. Hence $\overline{D_{2c_0}(\mathcal{C})}$ is reconstructible.

Finally, if the sum of two formal sums is reconstructible and one summand is
reconstructible, then so is the other. Applying this to Eq.~\eqref{eq:acover}
with $k=2c_0-1,2c_0-3,\ldots,1$, we obtain in turn the reconstructibility of
$\overline{D_{2c_0-2}(\mathcal{C})}$, $\overline{D_{2c_0-4}(\mathcal{C})}$, \ldots,
$\overline{D_0(\mathcal{C})}=\mathcal{C}$.
\end{proof}

\begin{theorem}\label{th:nuovo}
Let $(\Gamma,\sigma)$ be a connected balanced signed graph on $n\neq4$ vertices
whose balance partition $\{A,B\}$ satisfies $\bigl||A|-|B|\bigr|>2$. Then
$(\Gamma,\sigma)$ is reconstructible. Equivalently, a connected bipartitioned graph
$(\Gamma,\{A,B\})$ on $n\neq 4$ vertices can be non-reconstructible only if
$\{|A|,|B|\}=\{n/2,n/2\}$ or $\{|A|,|B|\}=\{n/2-1,n/2+1\}$.
\end{theorem}
\begin{proof}
By Proposition~\ref{prop:bip}, it suffices to show that $(\Gamma,\pi)$, with 
$\pi=\{A,B\}$, is reconstructible as a bipartitioned graph. Let $c_\pi$ be the
$2$-vertex-coloring of $\Gamma$ that colors in black the larger of $A$ and $B$.
The
$2$-vertex-colored deck of $(\Gamma,c_\pi)$ is determined by the bipartitioned
deck of $(\Gamma,\pi)$. In fact, since the parts of $\pi$ differ in size by more than $2$, in every card
$(\Gamma,\pi_v)$ the two parts still have different sizes, so that coloring in black
the larger part of $\pi_v$ gives exactly the card $(\Gamma,(c_\pi)_v)$.  Now let $(\Gamma,\pi')$ have the same bipartitioned deck as
$(\Gamma,\pi)$. By Lemma~\ref{prop:bic} the parts of $\pi'$ have the same sizes as
those of $\pi$, so $c_{\pi'}$ is defined in the same way and $(\Gamma,c_{\pi'})$
has the same $2$-vertex-colored deck as $(\Gamma,c_\pi)$. By
Proposition~\ref{teo:2}, $(\Gamma,c_{\pi'})\cong(\Gamma,c_\pi)$, hence
$(\Gamma,\pi')\cong(\Gamma,\pi)$.
\end{proof}

\begin{remark}\label{rem:aut}
In the case $\{|A|,|B|\}=\{n/2-1,n/2+1\}$, a natural idea is to use only the
$n/2-1$ cards obtained by switching a vertex of the smaller part $A$: these form
the deck of the subset $A\subseteq V_\Gamma$ in the sense of
\cite[Section~3.2]{ACKR}, and \cite[Corollary~2.4]{ACKR} gives the
reconstructibility of $(\Gamma,\pi)$ when $|\mathrm{Aut}(\Gamma)|<2^{|A|-1}=2^{n/2-2}$.
This is however precisely the bound under which \cite[Theorem~5.18]{Tan} gives
the reconstructibility of every signed graph on $\Gamma$, balanced or not.
\end{remark}
\subsection{Net-degrees in balanced signed graphs}\label{sec:netbal}
We show that, for connected balanced signed graphs, the case of positive
net-degrees of Theorem~\ref{cor:p1} extends to non-negative net-degrees
(Theorem~\ref{theo:final}). Two elementary facts about switch partners $v,w$ are
used repeatedly: how the net-degrees of the other vertices change from $\sigma$
to $\sigma_{vw}$ (Lemma~\ref{lem:swap}, valid for any two vertices), and which
part of the balance partition $w$ belongs to (Lemma~\ref{cor:n21}).

\begin{lemma}\label{lem:swap}
Let $v,w$ be distinct vertices of a signed graph $(\Gamma,\sigma)$. Then, for every
$x\notin\{v,w\}$,
$$
d^\pm_{\sigma_{vw}}(x)=d^\pm_\sigma(x)-2\bigl(\sigma(v,x)+\sigma(w,x)\bigr).
$$
\end{lemma}
\begin{proof}
After switching at $v$ and $w$, the only edges incident to $x$ that have changed sign are
$\{v,x\}$ and $\{w,x\}$, if present.
\end{proof}

\begin{lemma}\label{cor:n21}
Let $(\Gamma,\sigma)$ be a connected balanced signed graph on $n\neq4$ vertices and
let $v,w$ be switch partners. If $|\pi_\sigma(v)|\neq n/2+1$, then
$w\notin\pi_\sigma(v)$; in particular, if $v\sim w$ then $\sigma(v,w)=-1$.
\end{lemma}
\begin{proof}
If $w\in\pi_\sigma(v)=A$, the balance partition of $\sigma_{vw}$ has parts of sizes
$\{|A|-2,|B|+2\}$, which must equal $\{|A|,|B|\}$ by Lemma~\ref{prop:bic}; this
forces $|A|=n/2+1$.
\end{proof}

The following lemma, the technical core of the paper, handles net-degree $0$.

\begin{lemma}\label{lem:tec}
Let $(\Gamma,\sigma)$ be a connected balanced signed graph on $n\neq4$ vertices.
If $d^\pm_\sigma(v)\in\{0,2\}$ for every $v\in V_\Gamma$, then $(\Gamma,\sigma)$ is
reconstructible.
\end{lemma}
\begin{proof}
By Proposition~\ref{prop:4} we may assume $n\equiv0\pmod 4$, so $n\ge8$. Suppose that
$(\Gamma,\sigma)$ is not reconstructible. Then every vertex has a switch partner
(Proposition~\ref{lem:swp}), both net-degrees $0$ and $2$ occur
(Theorem~\ref{coro:rego}), and, since $\sigma_{vw}$ has the same deck as $\sigma$
for any switch partners $v,w$, all net-degrees of $\sigma_{vw}$ lie in $\{0,2\}$
by Proposition~\ref{prop:basic}(iii); we use this fact and Lemma~\ref{lem:swap}
throughout without further mention.

Let $d^\pm_\sigma(v)=2$ and let $w$ be a switch
partner of $v$. Proposition~\ref{sw.par.form} gives $d^\pm_\sigma(w)=2\,\sigma(v,w)-2$, thus $d^\pm_\sigma(w)=0$ and
$\sigma(v,w)=1$, so $w\in\pi_\sigma(v)$ and $|\pi_\sigma(v)|=n/2+1$ by
Lemma~\ref{cor:n21}. The same argument applies to every signature $\tau$ with the
same deck as $\sigma$ (such a $\tau$ is not reconstructible and
$\deg^\pm(\Gamma,\tau)\subseteq\{0,2\}$), hence:
\begin{equation}\label{eq:claim0}
\text{if $\tau$ has the same deck as $\sigma$ and $d^\pm_\tau(v)=2$, then } |\pi_\tau(v)|=n/2+1.
\end{equation}
Let $A$ be the part of $\pi_\sigma$ of size $n/2+1$ and $B$ the other one, of size
$n/2-1$. All vertices of net-degree $2$ lie in $A$; write $A=A_2\sqcup A_0$, where $A_2$ is
the set of vertices of net-degree $2$ and $A_0=A\setminus A_2$ consists of vertices
of net-degree $0$, and set $n_2=|A_2|$. 
Edges with both endpoints in $A$ or both in $B$
are positive, edges with one endpoint in $A$ and one in $B$ are negative. Note that
for $v,w\in A$ the balance partition of $\sigma_{vw}$ is
$\{A\setminus\{v,w\},\,B\cup\{v,w\}\}$, while for $u\in B$, $w\in A$ it is
$\{(A\setminus\{w\})\cup\{u\},\,(B\setminus\{u\})\cup\{w\}\}$.

\smallskip\noindent Claim \textbf{(1)} \emph{If $v\in A_2$ and $w\in A_0$ are switch partners,
they have no common neighbour.} A common neighbour $y$ has $\sigma(v,y)=\sigma(w,y)$
by balance, since $\sigma(v,w)=1$. Thus by Lemma~\ref{lem:swap} its net-degree in $\sigma_{vw}$ would be
$d^\pm_\sigma(y)\pm4\notin\{0,2\}$.

\smallskip\noindent Claim \textbf{(2)} \emph{Every $v\in A_2$ has exactly one neighbour in
$A_0$, namely its switch partner $w$, which is therefore unique; and $w$ has no
neighbour in $A_0$.} If $w'\in A_0\setminus\{w\}$ were adjacent to $v$, then
$w'\nsim w$ by (1), and by Lemma~\ref{lem:swap} $w'$ would have net-degree $-2$
in $\sigma_{vw}$; the same if $w'$ were adjacent to $w$.

\smallskip\noindent Claim \textbf{(3)} \emph{Let $v_1,v_2\in A_2$ have switch partners
$w_1,w_2$. Then $w_1\neq w_2$ if and only if $v_1\sim v_2$.} If $w_1\neq w_2$ and
$v_1\nsim v_2$, then $v_2\nsim w_1$ by (2), so $v_2$ keeps net-degree $2$ in
$\sigma_{v_1w_1}$, while by Proposition~\ref{prop:bip} it lies in the part
$A\setminus\{v_1,w_1\}$ of $\pi_{\sigma_{v_1w_1}}$, which has size $n/2-1$, against
Eq.~\eqref{eq:claim0}. If $w_1=w_2=w$ and $v_1\sim v_2$, then $v_2$ is a common
neighbour of the switch partners $v_1$ and $w$, against (1).

\smallskip\noindent Claim \textbf{(4)} \emph{If $u\in B$ and $w$ is a switch partner of $u$,
then $w\in A_0$ and $u\nsim w$.} By Lemma~\ref{cor:n21}, $w\notin B$;
Proposition~\ref{sw.par.form} gives $d^\pm_\sigma(w)=2\sigma(u,w)\in\{0,-2\}$,
so $d^\pm_\sigma(w)=0$ and $u\nsim w$.

\smallskip\noindent Claim \textbf{(5)} \emph{If $u\in B$ and $w\in A_0$ are switch partners,
then $\sigma(u,x)=1$ if and only if $\sigma(w,x)=-1$, for every $x\in V_\Gamma$.} If
$\sigma(u,x)=1$ then $x\in B$ loses $2$ from $u$, so it must gain $2$ from $w$. If
$\sigma(w,x)=-1$ then $x\in B$; if $x\nsim u$, then $x$ has net-degree $2$ in
$\sigma_{uw}$ while lying in the part $(B\setminus\{u\})\cup\{w\}$ of size $n/2-1$,
against Eq.~\eqref{eq:claim0}.

\smallskip\noindent Claim \textbf{(6)} \emph{If $u\in B$ and $w\in A_0$ are switch partners,
$x\in A_2$ and $\sigma(u,x)=-1$, then $\sigma(w,x)=1$.} Otherwise $x$ would have
net-degree $4$ in $\sigma_{uw}$.

\smallskip\noindent Claim \textbf{(7)} \emph{Every $w\in A_0$ has a neighbour in $A_2$;
consequently, by (2), all positive neighbours of $w$ lie in $A_2$.} Suppose $w$ has
no neighbour in $A_2$, and let $u$ be a switch partner of $w$. By
Proposition~\ref{sw.par.form}, $d^\pm_\sigma(u)=2\sigma(w,u)$, and $u\notin A_2$
(otherwise $\sigma(w,u)=1$ and $u$ would be a neighbour of $w$ in $A_2$); hence
$d^\pm_\sigma(u)=0$ and $u\nsim w$. Since $w$ is not isolated and
$d^\pm_\sigma(w)=0$, $w$ has a positive neighbour $x$, necessarily in $A_0$; $x$
loses $2$ in $\sigma_w$ and must gain $2$ from $u$, so $\sigma(u,x)=-1$ and $u\in B$.
Thus (5) and (6) apply to $u,w$. We claim that, since $w$ has no neighbour in $A_2$, $\sigma(u,x)=-1$ if and only if
$\sigma(w,x)=1$, for every $x\in V_\Gamma$. If $\sigma(w,x)=1$ then $x\in A_0$ and,
as just seen, $\sigma(u,x)=-1$. Conversely let $\sigma(u,x)=-1$, so $x\in A$: if
$x\in A_2$, (6) gives $\sigma(w,x)=1$; if $x\in A_0$ and $w\nsim x$, then $x$ gets
net-degree $2$ in $\sigma_{wu}$, while no vertex of $A_2$ changes its net-degree
(neither $w$ nor, by (6), $u$ has neighbours in $A_2$), so the number of vertices of
net-degree $2$ increases, against Proposition~\ref{prop:basic}(iii). Together with
(5), the claim says that in $(\Gamma,\sigma_w)$ the vertices $w$ and $u$ have the
same positive and the same negative neighbours outside $\{u,w\}$, so
Lemma~\ref{lemma:intorno} gives
$(\Gamma,\sigma_{wu})\cong(\Gamma,\sigma_{ww})=(\Gamma,\sigma)$, contradicting the
fact that $w,u$ are switch partners.

\smallskip
For $w\in A_0$ let $A_2^{(w)}$ be the set of positive neighbours of $w$, which lies
in $A_2$ by (7). By (2), (3) and (7), $A_2$ is the disjoint union of the sets
$A_2^{(w)}$, $w\in A_0$, each $v\in A_2$ having as unique switch partner the $w$ with
$v\in A_2^{(w)}$; the graph induced on $A_2$ is complete multipartite with parts
$A_2^{(w)}$; and
\begin{equation}\label{eq:PA}
n_2=\sum_{w\in A_0}d^+_\sigma(w).
\end{equation}

\smallskip\noindent Claim \textbf{(8)} \emph{$|A_0|\ge2$.} Suppose $A_0=\{w\}$. By (4),
$w$ is the switch partner of every vertex of $B$. Since $B\neq\emptyset$, we can
pick $u_1\in B$; it has net-degree $0$ and is not isolated, so it has a positive
neighbour $u_2$, which lies in $B$. By (5) applied to $u_1$ and $w$,
$\sigma(w,u_2)=-1$, contradicting (4) for $u_2$.

\smallskip\noindent Claim \textbf{(9)} \emph{$d^+_\sigma(w)\ge2$ for every $w\in A_0$.}
By (7), every term in Eq.~\eqref{eq:PA} is at least $1$, so $|A_0|\le n_2$ and
$n/2+1=|A|=n_2+|A_0|\le2n_2$; hence $n_2\ge3$, because $n\ge8$. Now suppose that $d^+_\sigma(w)=1$, say $A_2^{(w)}=\{v\}$, so that $w$ is
the switch partner of $v$ by (2). By (3) every vertex of
$A_2\setminus\{v\}$ is adjacent to $v$; together with $w$, which by (2) is the only
neighbour of $v$ in $A_0$, this gives $d^+_\sigma(v)=n_2$
 and $v$ has
$n_2-2\ge1$ negative neighbours $u_1,\dots,u_{n_2-2}$, all in $B$. By (4), the
switch partner of each $u_i$ lies in $A_0$, and by (6) it is positively adjacent
to $v$; it is therefore $w$. Let $x$ be the unique negative neighbour of
$w$, which lies in $B$. By (5), the positive neighbours of $u_i$ are exactly the
negative neighbours of $w$, so $x$ is the unique positive neighbour of $u_i$;
since $d^\pm_\sigma(u_i)=0$, $v$ is the unique negative neighbour of $u_i$.
Finally, let $w''$ be a switch partner of $x$: by (4), $w''\in A_0$, and by (5),
since $\sigma(x,u_i)=1$, we get $\sigma(w'',u_i)=-1$, so $w''=v\notin A_0$, a
contradiction.

\smallskip
From (8), (9) and Eq.~\eqref{eq:PA} we get $2\le d^+_\sigma(w)\le n_2-2$ for every
$w\in A_0$.

\smallskip\noindent Claim \textbf{(10)} \emph{$|A_0|=2$.} Fix $w\in A_0$, let
$k=d^+_\sigma(w)=d^-_\sigma(w)$, $A_2^{(w)}=\{v_1,\dots,v_k\}$, and let
$x_1,\dots,x_k\in B$ be the negative neighbours of $w$. By (3), the neighbours of
$v_1$ in $A_2$ are exactly the vertices of $A_2\setminus A_2^{(w)}$, so
$d^+_\sigma(v_1)=n_2-k+1$ and $v_1$ has $\ell=n_2-k-1\ge1$ negative neighbours
$u_1,\dots,u_\ell\in B$, whose switch partner is $w$ by (4) and (6). By (5),
$\sigma(u_i,x_j)=1$ for all $i,j$. Let $w'$ be a switch partner of $x_1$: by (4),
$w'\in A_0$ and $w'\neq w$, and (5) applied to $x_1$ gives $\sigma(w',u_i)=-1$ for
all $i$, so $d^+_\sigma(w')\ge\ell=n_2-k-1$. Thus $d^+_\sigma(w)+d^+_\sigma(w')\ge n_2-1$; since the sum in Eq.~\eqref{eq:PA}
equals $n_2$ and every term is at least $2$ by (9), a third vertex in $A_0$ would
make the sum exceed $n_2$. Hence $A_0=\{w,w'\}$ and $d^+_\sigma(w')=n_2-k$.

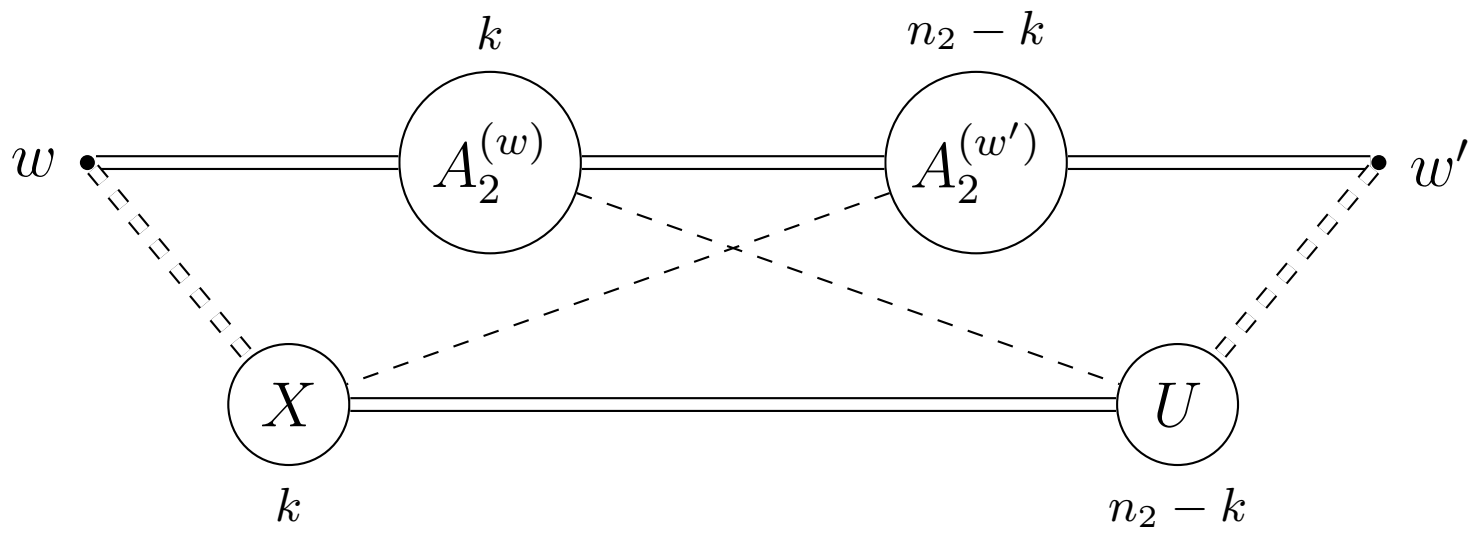
\begin{figure}
\centering
\begin{tikzpicture}[>=stealth, node distance=2.8cm and 2cm,scale=0.8]

  \node[fill=black, circle, inner sep=1pt, label=left:$w$] (w) {};
  \node[draw, circle, minimum size=1.2cm, label=center:$A_2^{(w)}$,
        label={[font=\scriptsize]90:$k$}] (A2w) [right=of w] {};
  \node[draw, circle, minimum size=1.2cm, label=center:$A_2^{(w')}$,
        label={[font=\scriptsize]90:$n_2-k$}] (A2wp) [right=of A2w] {};
  \node[fill=black, circle, inner sep=1pt, label=right:$w'$] (wp) [right=of A2wp] {};

  \node[draw, circle, minimum size=0.8cm, label=center:$X$,
        label={[font=\scriptsize]-90:$k$}] (X) at ($(w)!0.5!(A2w)+(0,-2)$) {};
  \node[draw, circle, minimum size=0.8cm, label=center:$U$,
        label={[font=\scriptsize]-90:$n_2-k$}] (U) at ($(A2wp)!0.5!(wp)+(0,-2)$) {};

  \draw[double, double distance=2pt] (w) -- (A2w);
  \draw[double, double distance=2pt] (A2w) -- (A2wp);
  \draw[double, double distance=2pt] (A2wp) -- (wp);
  \draw[double, double distance=2pt] (X) -- (U);
  \draw[dashed, double, double distance=2pt] (w) -- (X);
  \draw[dashed, double, double distance=2pt] (wp) -- (U);

  \draw[dashed] (A2w) -- (U);
  \draw[dashed] (A2wp) -- (X);

\end{tikzpicture}
\caption{The structure of $(\Gamma,\sigma)$ in the proof of Lemma~\ref{lem:tec},
once $A_0=\{w,w'\}$ has been established. A double line means that all edges
between the two ends are present, a single line that only some of them are, all
with the sign indicated by the line.}\label{fig:proof}
\end{figure}

 Let $A_0=\{w,w'\}$, with $k=d^+_\sigma(w)$ and
$d^+_\sigma(w')=n_2-k$ as in (10). Now $|A|=n_2+2$, so $n_2=n/2-1=|B|$, which is
odd. Let $X\subseteq B$ be the set of negative neighbours of $w$ and $U\subseteq B$
that of $w'$, of sizes $k$ and $n_2-k$. A vertex of $B$ is not adjacent to its switch
partner (4), so $X\cap U=\emptyset$ and $B=X\sqcup U$, every $x\in X$ having partner
$w'$ and every $u\in U$ having partner $w$; the situation is depicted in
Fig.~\ref{fig:proof}. By (5), the positive neighbours of $x\in X$ are exactly the
vertices of $U$, and those of $u\in U$ are exactly the vertices of $X$: there are no
edges inside $X$ or inside $U$, and every vertex of $U$ is adjacent to every vertex
of $X$. By (1), no vertex of $A_2^{(w')}$ is adjacent to a vertex of $U$ and no
vertex of $A_2^{(w)}$ is adjacent to a vertex of $X$. By (2) and (3), the positive
neighbours of $v\in A_2^{(w')}$ are exactly the vertices of $A_2^{(w)}\cup\{w'\}$, so
$d^+_\sigma(v)=k+1$, $d^-_\sigma(v)=k-1$, and all its negative neighbours lie in $X$.
Count the negative edges between $X$ and $A_2^{(w')}$: each $x\in X$ has
$d^-_\sigma(x)=d^+_\sigma(x)=n_2-k$ negative neighbours, namely $w$ and $n_2-k-1$
vertices of $A_2^{(w')}$, giving $k(n_2-k-1)$ edges; counting from $A_2^{(w')}$ gives
$(n_2-k)(k-1)$ edges. Equality forces $n_2=2k$, which is impossible since $n_2$ is
odd.
\end{proof}

\begin{theorem}\label{theo:final}
Let $(\Gamma,\sigma)$ be a connected balanced signed graph on $n\neq4$ vertices. If
all net-degrees are non-negative, then $(\Gamma,\sigma)$ is reconstructible.
\end{theorem}
\begin{proof}
Let $d_1>\dots>d_m$ be the distinct net-degrees. By Theorem~\ref{cor:p1} we may
assume $d_m=0$, and by Proposition~\ref{cor:d1dn} that $d_1\le2$. Suppose that
$(\Gamma,\sigma)$ is not reconstructible; by Lemma~\ref{lem:tec} some vertex $s$
has net-degree $1$.
A switch partner $s'$ of $s$ satisfies $1+d^\pm_\sigma(s')=2\sigma(s,s')$ by
Proposition~\ref{sw.par.form}, and since all net-degrees are non-negative this
forces $d^\pm_\sigma(s')=1$ and $\sigma(s,s')=1$; so $s'\in\pi_\sigma(s)$ and
$|\pi_\sigma(s)|=n/2+1$ by Lemma~\ref{cor:n21}. The same holds for every signature
with the same deck as $\sigma$.

\emph{Case 1: some $v$ has net-degree $2$.} Its switch partner $w$ has net-degree $0$
and $\sigma(v,w)=1$ by Proposition~\ref{sw.par.form} and since all net-degrees are non-negative, so $|\pi_\sigma(v)|=n/2+1$ by
Lemma~\ref{cor:n21}, and $s,v,w$ all lie in the part $A$ of size $n/2+1$. Let
$\tau=\sigma_{vw}$, which has the same deck as $\sigma$. Switching preserves
net-degrees modulo $2$, so $d^\pm_\tau(s)$ is odd, and since
$\deg^\pm(\Gamma,\tau)=\deg^\pm(\Gamma,\sigma)\subseteq\{0,1,2\}$ by
Proposition~\ref{prop:basic}(iii), $d^\pm_\tau(s)=1$. But by
Proposition~\ref{prop:bip} $s$ lies in the part $A\setminus\{v,w\}$ of $\pi_\tau$,
which has size $n/2-1$, contradicting what we proved for $\tau$.

\emph{Case 2: all net-degrees are $0$ or $1$.} As $\Gamma$ is connected and both values
occur, there are adjacent vertices $s,w$ with $d^\pm_\sigma(s)=1$ and
$d^\pm_\sigma(w)=0$. Let $s'$ be a switch partner of $s$. In $\sigma_{ss'}$ the vertex
$w$ has an even net-degree in $\{0,1\}$, hence $0$, so by Lemma~\ref{lem:swap}
$\sigma(s,w)+\sigma(s',w)=0$: $w$ is adjacent to both $s$ and $s'$,
with $\sigma(s,w)=-\sigma(s',w)$. This is impossible, since $s$ and $s'$ lie in the
same part and $\sigma$ is balanced.
\end{proof}

\begin{remark}
In the language of Section~\ref{sec:bip}, Theorem~\ref{theo:final} states that a
connected bipartitioned graph on $n\neq4$ vertices is reconstructible whenever
every vertex has at least as many neighbours in its own part as in the other part.
By Proposition~\ref{prop:comple}, the same holds when every vertex has at least as
many non-neighbours in its own part as in the other part, provided $\Gamma^c$ is
connected.
\end{remark}

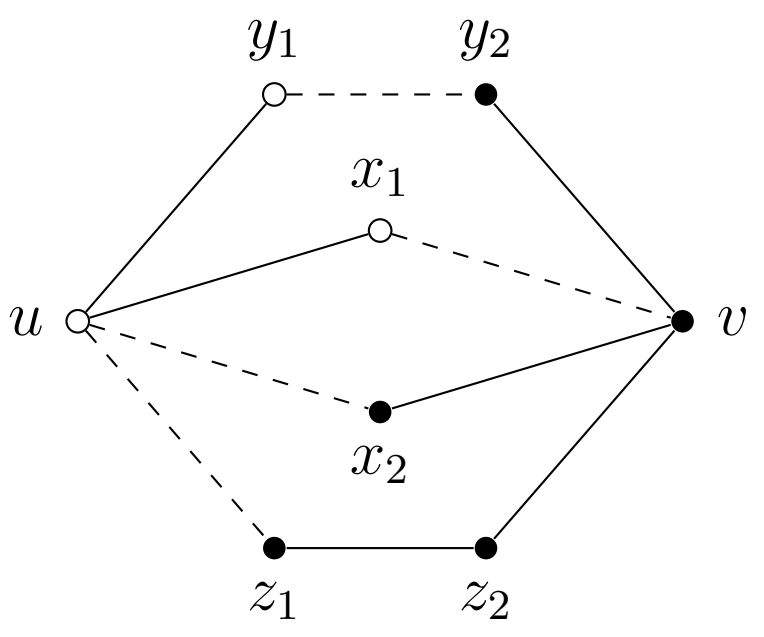
\begin{figure}
\centering
\begin{tikzpicture}[scale=1]
  \node[draw, circle, inner sep=1.5pt, label=left:$u$]        (u)  at (-2,0) {};
  \node[fill=black, circle, inner sep=1.5pt, label=right:$v$] (v)  at (2,0) {};
  \node[draw, circle, inner sep=1.5pt, label=above:$y_1$]       (y1) at (-0.7,1.5) {};
  \node[fill=black, circle, inner sep=1.5pt, label=above:$y_2$] (y2) at (0.7,1.5) {};
  \node[draw, circle, inner sep=1.5pt, label=above:$x_1$]       (x1) at (0,0.6) {};
  \node[fill=black, circle, inner sep=1.5pt, label=below:$x_2$] (x2) at (0,-0.6) {};
  \node[fill=black, circle, inner sep=1.5pt, label=below:$z_1$] (z1) at (-0.7,-1.5) {};
  \node[fill=black, circle, inner sep=1.5pt, label=below:$z_2$] (z2) at (0.7,-1.5) {};
  \draw (u) -- (x1);  \draw[dashed] (x1) -- (v);
  \draw[dashed] (u) -- (x2);  \draw (x2) -- (v);
  \draw (u) -- (y1);  \draw[dashed] (y1) -- (y2);  \draw (y2) -- (v);
  \draw[dashed] (u) -- (z1);  \draw (z1) -- (z2) -- (v);
\end{tikzpicture}
\caption{The signed graph of Example~\ref{exa:theta}: filled vertices form $A$,
hollow vertices form $B$.}\label{fig:theta}
\end{figure}

\begin{example}\label{exa:theta}
Let $\Gamma$ be the graph on $8$ vertices consisting of two vertices $u,v$ joined
by two paths of length $2$, $u\,x_1\,v$ and $u\,x_2\,v$, and two paths of length
$3$, $u\,y_1y_2\,v$ and $u\,z_1z_2\,v$, and let $\sigma$ be the balanced signature
with balance partition
$$
A=\{v,x_2,y_2,z_1,z_2\},\qquad B=\{u,x_1,y_1\},
$$
see Fig.~\ref{fig:theta}. The parts have sizes $5$ and $3$, so
Theorem~\ref{th:nuovo} does not apply. The vertex $v$ has three neighbours in $A$
and one in $B$, and $z_2$ has both its neighbours in $A$, so they have net-degree
$2$; every other vertex has net-degree $0$.
Thus  $(\Gamma,\sigma)$ is
reconstructible by Theorem~\ref{theo:final}. Finally, for each of the subsets
$$
\{x_1,x_2,y_1,z_1\},\ \{x_1,x_2,z_1,z_2\},\ \{x_1,x_2,y_2,z_1\},\
\{u,v,y_1,z_2\},\ \{u,v,y_1,y_2\},\ \{u,v,y_2,z_2\}
$$
one checks that $(\Gamma,\sigma_W)\cong(\Gamma,\sigma)$, and each vertex lies in
exactly three of them. Together with $\emptyset$ and $V_\Gamma$, these subsets
 meet the lower bounds required by \cite[Theorem~5.17]{Tan} for a
non-reconstructible signed graph on $8$ vertices.  Therefore  \cite[Theorem~5.17]{Tan},
together with its consequence \cite[Theorem~5.18]{Tan} discussed in
Remark~\ref{rem:aut}, gives no information on $(\Gamma,\sigma)$.
\end{example}

\begin{corollary}\label{coro:bip}
Let $(\Gamma,\sigma)$ be a connected balanced signed graph on $n\neq4$ vertices
with $\Gamma$ bipartite. If no vertex has positive
net-degree, or no vertex has negative net-degree, then $(\Gamma,\sigma)$ is
reconstructible.
\end{corollary}
\begin{proof}
If no net-degree is negative, apply Theorem~\ref{theo:final}. Otherwise consider
$(\Gamma,-\sigma)$, which is balanced because $\Gamma$ is bipartite
\cite[Corollary~I.3]{Zaslavsky}, has no negative net-degree, and is reconstructible
if and only if $(\Gamma,\sigma)$ is.
\end{proof}
Signed trees are balanced and bipartite, so Corollary~\ref{coro:bip} applies to
them in particular.

\section*{Acknowledgements and declarations}
The authors acknowledge the support of the INdAM-GNSAGA project \lq\lq Between groups and graphs\rq\rq, CUP E53C25002010001. This research was partially supported by the Engineering Department of Cusano University of Rome.
A. Donno was partially supported by the Research Project \lq\lq Operator algebras and harmonic analysis, with
applications to noncommutative geometry, quantum field theory, representation theory and combinatorics\rq\rq funded by Sapienza University of Rome.
\\The authors declare that there are no competing interests that are directly or indirectly related to the work submitted for publication.
\\Data sharing is not applicable to this article as no datasets were generated or analyzed during the current study.


\begin{thebibliography}{99}

\bibitem{ACKR} N. Alon, Y. Caro, I. Krasikov, Y. Roditty, Combinatorial reconstruction problems, \emph{J. Combin. Theory Ser. B} {\bf 47} (1989), no. 2, 153--161.

\bibitem{bo} J. A. Bondy, Beautiful conjectures in graph theory, \emph{European J. Combin.} {\bf 37} (2014), 4--23.

\bibitem{BondyHemminger} J. A. Bondy, R. L. Hemminger, Graph reconstruction--a survey, \emph{J. Graph Theory} {\bf 1} (1977), no. 3, 227--268.

\bibitem{BondyMercier} J. A. Bondy, F. Mercier, Switching reconstruction of digraphs, \emph{J. Graph Theory} {\bf 67} (2011), no. 4, 332--348.

\bibitem{CDSdm} M. Cavaleri, A. Donno, S. Spessato, On the complement of a signed graph, \emph{Discrete Math.} {\bf 348} (2025), no. 6, Paper No.~114433, 21 pp.

\bibitem{Ell} M. N. Ellingham, Vertex-switching reconstruction and folded cubes, \emph{J. Combin. Theory Ser. B} {\bf 66} (1996), no. 2, 361--364.

\bibitem{EllinghamRoyle} M. N. Ellingham, G. F. Royle, Vertex-switching reconstruction of subgraph numbers and triangle-free graphs, \emph{J. Combin. Theory Ser. B} {\bf 54} (1992), no. 2, 167--177.

\bibitem{Harary} F. Harary, On the notion of balance of a signed graph, \emph{Michigan Math. J.} {\bf 2} (1953--1954), 143--146; addendum, \emph{Michigan Math. J.} {\bf 2}, preceding p.~1.

\bibitem{Harary1985} F. Harary, The reconstruction conjecture for balanced signed graphs, in \emph{Cycles in graphs} (Burnaby, B.C., 1982), pp. 439--442, North-Holland Math. Stud., 115, Ann. Discrete Math., 27, North-Holland, Amsterdam, 1985.

\bibitem{kra88} I. Krasikov, A note on the vertex-switching reconstruction, \emph{Internat. J. Math. Math. Sci.} {\bf 11} (1988), no. 4, 825--827.

\bibitem{kra94} I. Krasikov, Applications of balance equations to vertex switching reconstruction, \emph{J. Graph Theory} {\bf 18} (1994), no. 3, 217--225.

\bibitem{kra96} I. Krasikov, Degree conditions for vertex switching reconstruction, \emph{Discrete Math.} {\bf 160} (1996), no. 1-3, 273--278.

\bibitem{KrasikovRoditty1987} I. Krasikov, Y. Roditty, Balance equations for reconstruction problems, \emph{Arch. Math. (Basel)} {\bf 48} (1987), no. 5, 458--464.

\bibitem{KrasikovRoditty1992} I. Krasikov, Y. Roditty, Switching reconstruction and Diophantine equations, \emph{J. Combin. Theory Ser. B} {\bf 54} (1992), no. 2, 189--195.

\bibitem{KrasikovRoditty1994} I. Krasikov, Y. Roditty, More on vertex-switching reconstruction, \emph{J. Combin. Theory Ser. B} {\bf 60} (1994), no. 1, 40--55.

\bibitem{McKaySchweitzer} B. D. McKay, P. Schweitzer, Switching reconstruction of digraphs, \emph{J. Graph Theory} {\bf 76} (2014), no. 4, 279--296.

\bibitem{Mc19} J. C. McLeod, B. D. McKay, B. Faller, The $r$-switching-stable graphs, \emph{Discrete Appl. Math.} {\bf 266} (2019), 16--29.

\bibitem{Nie} P. Niesink, The vertex-switching reconstruction problem, M.Sc. thesis, University of Ottawa, 2010.

\bibitem{sei} J. J. Seidel, A survey of two-graphs, in \emph{Colloquio Internazionale sulle Teorie Combinatorie (Rome, 1973), Tomo I}, pp. 481--511, Atti dei Convegni Lincei, No. 17, Accad. Naz. Lincei, Rome, 1976.

\bibitem{Stanley} R. P. Stanley, Reconstruction from vertex-switching, \emph{J. Combin. Theory Ser. B} {\bf 38} (1985), no. 2, 132--138.

\bibitem{Tan} J. Tan, Graph reconstruction and structure, DPhil thesis, University of Oxford, 2023.

\bibitem{zasgraph} T. Zaslavsky, Characterizations of signed graphs, \emph{J. Graph Theory} {\bf 5} (1981), no. 4, 401--406.

\bibitem{ZaslavskyBiblio} T. Zaslavsky, A mathematical bibliography of signed and gain graphs and allied areas, \emph{Electron. J. Combin.} {\bf 5} (1998), Dynamic Surveys 8, 124 pp.

\bibitem{Zaslavsky} T. Zaslavsky, Matrices in the theory of signed simple graphs, in \emph{Advances in discrete mathematics and applications: Mysore, 2008}, pp. 207--229, Ramanujan Math. Soc. Lect. Notes Ser. {\bf 13}, Ramanujan Mathematical Society, Mysore, 2010.

\end{thebibliography}
\end{document}